\documentclass[english]{article}
\usepackage[T1]{fontenc}
\usepackage[latin9]{inputenc}
\usepackage{amsmath}
\usepackage{amssymb}
\usepackage{graphicx}
\usepackage{babel}
\usepackage{xcolor}

\title{Clustering of Series via Dynamic Mode Decomposition and the Matrix Pencil Method}
\author{Leonid Pogorelyuk and Clarence W. Rowley}

\begin{document}

\maketitle

\begin{abstract}
In this paper, a new algorithm for extracting features from sequences of multidimensional observations is presented. The independently developed Dynamic Mode Decomposition and Matrix Pencil methods provide a least-squares model-based approach for estimating complex frequencies present in signals as well as their corresponding amplitudes. Unlike other feature extraction methods such as Fourier Transform or Autoregression which have to be computed for each sequence individually, the least-squares approach considers the whole dataset at once. It invokes order reduction methods to extract a small number of features best describing all given data, and indicate which frequencies correspond to which sequences. As an illustrative example, the new method is applied to regions of different grain orientation in a Transmission Electron Microscopy image.
\end{abstract}

\section{Introduction}

Clustering of series of data points is an unsupervised classification
task which consists of grouping sequences together based on some notion
of similarity between them. In the literature this task is known as ``clustering
of time series'' or ``functional data clustering'', although it
applies to any sequences, not just sequences in time \cite{liao2005clustering,jacques2014functional}.
These methods belong to the field of cluster analysis, which has an
even larger variety of algorithms and techniques \cite{kaufman2009finding}.

It is common to divide the various time series clustering algorithms
into several approaches \cite{liao2005clustering}: The raw-data-based
approaches compute similarity or distance metrics directly on the
time series and then apply more general clustering algorithms (e.g.
k-means \cite{macqueen1967some}) to achieve their goal. Feature-based
approaches first compute features using methods such as the Fast Fourier Transform (FFT)
\cite{welch1967use} or Principal Component Analysis (PCA) \cite{hotelling1933analysis}
on the time series and then apply standard clustering techniques. Finally, model-based
approaches assume some model for the process producing the data series
and compute its parameters.

A common model-based approach assumes some linear relation between
elements of the series, their precursors and process noise \cite{piccolo1990distance,xiong2002mixtures}.
In that case, one fits autoregression coefficients to the series first
\cite{akaike1969fitting}, and then proceeds by computing a distance metric
between time series and clustering based on those coefficients. Recently,
\cite{surana2018koopman} has introduced another metric for time series, which assumes a linear relation between
elements within each series, and relies on a method called Dynamic Mode Decomposition
(DMD) \cite{schmid2010dynamic}.

DMD \cite{schmid2010dynamic} and its variants assume that the data
was generated by a linear dynamical system with measurement noise; the method then estimates the complex frequencies and magnitudes (called ``modes'')
corresponding to that system. For sequential data, DMD is often paired with
delay-embedding, a common technique in system identification \cite{juang1985eigensystem}.
It involves defining high dimensional ``delayed observables'' which
include several adjacent data points from the given sequence, arranging
the data in (generalized) Hankel matrices and then applying standard
DMD \cite{tu2013dynamic,arbabi2016ergodic,zhang2017online}.

Interestingly, Hankel matrices were employed in the literature of
frequency estimation long before the appearance of DMD (for a review
see \cite{hokanson2013numerically}). Modern algorithms for frequency
estimation include the State Space method \cite{kung1983state}, ESPRIT
\cite{roy1986estimation} and the the Matrix Pencil method \cite{hua1988matrix}.
Similarly to DMD, they all compute the Singular Value Decomposition
(SVD) of the Hankel data matrices and, under certain conditions to
be discussed in this paper, give the same results.

Furthermore, the various methods---Matrix Pencil, ESPRIT, and DMD---each provide
a different perspective on the same arrangement of given data series.
Together, they offer an insight on designing a new set of features
for sequential data which will be the main subject of this paper.
Those features, computed on a set of series, are particularly useful for
clustering purposes.

In section~\ref{sec:DMD_MP_comparison}, we give an overview of DMD, Matrix Pencil methods,
and ESPRIT, emphasizing their similarities. Based
on these, an algorithm for extracting features from a set of sequences
in described in section~\ref{sec:Clustering}. Finally, numerical
examples for the use of those features for clustering are presented
in section~\ref{sec:Examples}.

\section{\label{sec:DMD_MP_comparison}Comparison Between DMD and Matrix Pencil
Methods}

Throughout our discussion, we consider data consisting of measurements
$y(t)$ of a signal $x(t)$ in the presence of noise $s(t)$,
where $t=0,1,\ldots,T$.  We also assume that the signal is a linear combination of $R$ damped sinusoids, so that it may be written
\begin{equation}
\begin{aligned}
  x(t)&=\underset{j=1}{\overset{R}{\sum}}v_{j}\lambda_{j}^{t}\\
  y(t)&=x(t)+s(t),
\end{aligned}
\label{eq:signal}
\end{equation}
where $x(t),y(t),s(t),v_{j}\in\mathbb{C}^{n}$, $\lambda_{j}\in\mathbb{C}$. In this section we discuss two families of methods
which estimate $R$ and $\lambda_{j},v_{j}$ based on the data $\left\{ y(t)\right\} _{t=0}^{T}$.

Earlier frequency estimation methods first approximate the complex frequencies
$\lambda_1,\ldots,\lambda_R$, and leave the task of finding the coefficients
$v_1,\ldots,v_R$ for a standard least squares procedure. Among these
are the ESPRIT \cite{roy1986estimation}, the Matrix Pencil \cite{hua1990matrix}
and State Space \cite{kung1983state} methods, which will be discussed
in sections \ref{sub:matrix_pencil}--\ref{sub:ESPRIT}.

The Dynamic Mode Decomposition (DMD) \cite{schmid2010dynamic} and
its variants estimate both the frequencies and the (vector) coefficients $v_1,\ldots,v_R$
(called ``modes'') and are frequently applied to non-sequential
data as well. In the context of sequential data, DMD is often used
with delayed observables \cite{tu2013dynamic,arbabi2016ergodic,zhang2017online},
which will be defined in sections \ref{sub:DMD}-\ref{sub:DMD_implementation}.

This section will focus on showing that DMD with delayed observables
is equivalent to the Matrix Pencil methods which historically preceded
it, and will consider a possible insight from the ESPRIT method to
the DMD community.

\subsection{\label{sub:DMD}DMD with Delayed Observables and Order Reduction}

The Dynamic Mode Decomposition (DMD) was originally
  proposed in the fluid mechanics community, as a method for identifying coherent structures in fluids flows. In DMD setting, each pair of
measurements $\left(y(t),y(t)\right)$ from (\ref{eq:signal})
is called a data snapshot, with the indices $t$ not necessarily
being ordered or sequential.  The general idea behind DMD is to find a linear map such that
\begin{equation}
y(t+1)\approx Ly(t),
\end{equation}
with $L\in\mathbb{C}^{n\times n}$. However, since $L$ has at most
rank $n$, this approximation is extremely poor when $n<R$ even if
the measurements are exact, i.e. $s(t)=0,\:\forall t$. To alleviate
this problem in the case when sequential data is available, it is
common to introduce delayed observables \cite{juang1985eigensystem}, as
\begin{equation}
z(t)=\begin{bmatrix}
y(t)\\
y(t+1)\\
\vdots\\
y(t+d)
\end{bmatrix}\in\mathbb{C}^{\left(d+1\right)n},\label{eq:delayed_observable}
\end{equation}
with $d>0$ delays and $\left(d+1\right)n\ge R$. When no noise is
present ($s(t)=0$), given enough data, one can fit an autoregressive
model \cite{akaike1969fitting}, such that
\begin{equation}
z(t)=Kz(t),\label{eq:delayed_map}
\end{equation}
holds exactly for some $K\in\mathbb{C}^{\left(d+1\right)n\times\left(d+1\right)n}$.
In this case, $R$ eigenvalues of $K$ are $\lambda_1,\ldots,\lambda_R$, and
the rest are zero.

If $R$ is not precisely known but $\left(d+1\right)n>R$, it is necessary
to get an estimate of the number of frequencies, $\bar{R}$, and project
the delayed observables $z(t)$ onto an $\bar{R}<\left(d+1\right)n$
dimensional space,
\begin{equation}
\tilde{z}(t)=Tz(t),
\end{equation}
where $\tilde{z}(t)$ are the reduced order observables and $T\in\mathbb{C}^{\bar{R}\times\left(d+1\right)n}$
has rank $\bar{R}$. In the DMD literature it is customary to compute the
reduced-order observables using
Proper Orthogonal Decomposition (POD)
\cite{lumley1970stochastic,schmid2010dynamic} (also known as Principal Component
Analysis (PCA) \cite{hotelling1933analysis,kung1983state}), in which $T$ is
computed using Singular Value
Decomposition (SVD), as described in the following subsection.

Consequently, considering the data pair $\left(\tilde{z}(t),\tilde{z}(t+1)\right)$,
the fitted linear map
\begin{equation}
\tilde{z}(t+1)=\tilde{K}\tilde{z}(t),
\end{equation}
with $\tilde{K}\in\mathbb{C}^{\bar{R}\times\bar{R}}$ will, assuming
$\bar{R}=R$ and $s(t)=0$, have exactly the desired eigenvalues
$\lambda_1,\ldots,\lambda_R$.

Besides estimating the frequency, one can decompose the signal into
a sum,
\begin{equation}
z(t)=\underset{j=1}{\overset{R}{\sum}}\begin{bmatrix}
v_{j}\\
\lambda_{j}v_{j}\\
\vdots\\
\lambda_{j}^{d}v_{j}
\end{bmatrix}\lambda_{j}^{t}+\begin{bmatrix}
s(t)\\
s(t+1)\\
\vdots\\
s(t+d)
\end{bmatrix},\label{eq:delayed_dynamics}
\end{equation}
of so called dynamic modes, $\begin{bmatrix}
v_{j}^{T} & ... & \lambda_{j}^{d}v_{j}^{T}\end{bmatrix}^{T}$, and a remainder, $\begin{bmatrix}
s(t)^{T} & \cdots & s(t+d)^{T}\end{bmatrix}^{T}$. Those DMD modes are eigenvectors of $K$ in (\ref{eq:delayed_map})
and a discussion on their scaling can be found in \cite{tu2013dynamic}.

We next present the full implementation details of the DMD method.

\subsection{\label{sub:DMD_implementation}DMD implementation}

Given $M$ pairs of snapshots $\left(z(t),z(t+1)\right)$, the data matrices are defined such that their columns consist of the delayed observables
in no particular order,
\begin{equation}
X=\begin{bmatrix}
  z(t_{0}) & \cdots & z(t_{M-1})\end{bmatrix},\qquad
Y=\begin{bmatrix}
z(t_{0}+1) & \cdots & z(t_{M-1}+1)
\end{bmatrix}
.\label{eq:data_matrices}
\end{equation}
For a sequential signal, one may choose $t_{i}=i$ and construct the data matrices,
\begin{equation}
X=\begin{bmatrix}
y(0) & y(1) & \cdots\\
y(1) & y(2)\\
\vdots & \vdots\\
y(d) & y(d+1) & \cdots
\end{bmatrix}
,\:Y=\begin{bmatrix}
y(1) & y(2) & \cdots\\
y(2) & y(3)\\
\vdots & \vdots\\
y(d+1) & y(d+2) & \cdots
\end{bmatrix}
,\label{eq:sequential_data_matrices}
\end{equation}
to be used in the following algorithm.

\emph{Step 1}: Compute the Singular Value Decomposition (SVD) 
\begin{equation}
X=U\Sigma V^{*},
\end{equation}
and choose the largest $\bar{R}$ singular values after which there
is a significant drop in their magnitude. The truncated SVD approximation
of the data is then $X\approx U_{\bar{R}}\Sigma_{\bar{R}}V_{\bar{R}}^{*}$
with $U_{\bar{R}}\in\mathbb{C}^{\left(d+1\right)n\times\bar{R}}$,
$\Sigma_{\bar{R}}\in\mathbb{C}^{\bar{R}\times\bar{R}}$, $V_{\bar{R}}\in\mathbb{C}^{M\times\bar{R}}$,
and $T=U_{R}^{*}$ is the order reducing transformation.

\emph{Step 2}:\emph{ }Defining the reduced order data matrices as
$\tilde{X}=U_{\bar{R}}^{*}X$, $\tilde{Y}=U_{\bar{R}}^{*}Y$ both
in $\mathbb{C}^{\bar{R}\times M}$, we have
\begin{equation}
\tilde{K}=\tilde{Y}\tilde{X}^{+}=U_{\bar{R}}^{*}YV_{\bar{R}}\Sigma_{\bar{R}}^{-1},\label{eq:DMD_matrix}
\end{equation}
where $\tilde{K}\in\mathbb{C}^{\bar{R}\times\bar{R}}$.

\emph{Step 3}: Perform an eigendecomposition of $\tilde{K}$
\begin{equation}
\tilde{K}=\tilde{V}\bar{\Lambda}\tilde{V}^{-1},\label{eq:eigendecomposition}
\end{equation}
where columns of $\tilde{V}\in\mathbb{C}^{\bar{R}\times\bar{R}}$ ($\tilde{V}^{-1}\in\mathbb{C}^{\bar{R}\times\bar{R}}$)
are the right (left) eigenvectors of $\tilde{K}$ and $\bar{\Lambda}=\mathrm{diag}\left\{ \bar{\lambda}_{1\bar{R}},...,\bar{\lambda}_{\bar{R}}\right\} $
consists of the approximated complex frequencies.

\emph{Step 4}: A rank $\bar{R}$ approximation of the full dynamics
in (\ref{eq:delayed_map}) is 
\begin{equation}
K\approx U_{\bar{R}}\tilde{K}U_{\bar{R}}^{*}.\label{eq:delayed_map_approx}
\end{equation}
Its right eigenvectors, the columns of $U_{\bar{R}}\tilde{V}$, are
called the DMD modes, while the rows of $\tilde{V}^{-1}U_{\bar{R}}^{*}$
are the left eigenvectors and are known as the adjoint DMD modes \cite{tu2013dynamic}.

\emph{Step }5: To find an appropriate scaling for the DMD modes we
first notice that for any $t,\tau$ it follows from (\ref{eq:delayed_map}),
(\ref{eq:eigendecomposition}) and (\ref{eq:delayed_map_approx})
that
\begin{equation}
z(t)\approx K^{t-\tau}z(\tau)\approx\underset{j=1}{\overset{\bar{R}}{\sum}}c_{\tau,j}w_{j}\bar{\lambda}_{j}^{t-\tau},
\end{equation}
where $w_{j}=\left(U_{\bar{R}}\tilde{V}\right)_{j}\in\mathbb{C}^{\left(d+1\right)n}$
is the $j$th DMD mode and $c_{\tau,j}=\left(\tilde{V}^{-1}U_{\bar{R}}^{*}z(\tau)\right)_{j}\in\mathbb{C}$
are scaling coefficients based on the $\tau$th snapshot.

In general, since $M>\bar{R}$, one has $c_{\tau_{1},j}\bar{\lambda}_{j}^{-\tau_{1}}\ne c_{\tau_{2},j}\bar{\lambda}_{j}^{-\tau_{2}}$
(unlike in the underconstrained case discussed in \cite{tu2013dynamic}),
hence we suggest a scaling for the $j$th DMD mode based on an average
of its scaling coefficients at all time snapshots,
\begin{equation}
c_{j}=\frac{1}{M}\underset{\tau\in\left\{ t_{1},...,t_{M}\right\} }{\sum}\left(\tilde{V}^{-1}U_{\bar{R}}^{*}z(\tau)\right)_{j}\bar{\lambda}_{j}^{-\tau}.
\end{equation}
This gives an approximation of the dynamics in delayed observables
(see (\ref{eq:delayed_dynamics})),
\begin{equation}
z(t)\approx\underset{j=1}{\overset{\bar{R}}{\sum}}c_{j}w_{j}\bar{\lambda}_{j}^{t},\label{eq:approx_delayed_dynamics}
\end{equation}
where $c_{j}w_{j}\in\mathbb{C}^{\left(d+1\right)n}$ are the scaled
dynamic modes of the delayed observable. The coefficients $v_{j}$
in (\ref{eq:signal}) are then estimated via
\begin{equation}
\bar{v}_{j}=\frac{1}{d+1}\underset{k=0}{\overset{d}{\sum}}\left(c_{j}w_{j}\right)_{kn:\left(k+1\right)n}\bar{\lambda}_{j}^{-k},
\end{equation}
where $\left(c_{j}w_{j}\right)_{kn:\left(k+1\right)n}\in\mathbb{C}^{n}$
are the $n$ consequent elements of $c_{j}w_{j}$ beginning at the
$kn$th elements.

\subsection{\label{sub:matrix_pencil}Matrix Pencil Methods}
The authors of the Matrix Pencil \cite{hua1988matrix} and the State Space \cite{kung1983state} methods originally considered the problem of retrieving parameters of sinusoidal processes (with frequencies close to one another) from noisy measurements. There, the data matrices are defined as 
\begin{equation}
X=\begin{bmatrix}
\cdots & y(1) & y(0)\\
 & y(2) & y(1)\\
 & \vdots & \vdots\\
\cdots & y(d+1) & y(d)
\end{bmatrix},\:Y=\begin{bmatrix}
\cdots & y(2) & y(1)\\
 & y(3) & y(2)\\
 & \vdots & \vdots\\
\cdots & y(d+2) & y(d+1)
\end{bmatrix},
\end{equation}
which correspond to a reversed ordering with respect to (\ref{eq:sequential_data_matrices}). The Matrix Pencil method proceeds by finding the generalized eigenvectors
of the matrix pencil $\left(X,Y\right)$, i.e. $\bar{\lambda}\in\mathbb{C}$
and $p\in\mathbb{C}^{M}$, $q\in\mathbb{C}^{\left(d+1\right)n}$ in
the row and column spaces of $X$ respectively, such that
\begin{equation}
\left(Y-\bar{\lambda}X\right)p=0,
\end{equation}
\begin{equation}
q^{*}\left(Y-\bar{\lambda}X\right)=0.
\end{equation}
Similarly to DMD, this is accomplished by computing the truncated
SVD of $X$, $X\approx U_{\bar{R}}\Sigma_{\bar{R}}V_{\bar{R}}^{*}$.
The generalized eigenvalues $\bar{\lambda}$ are estimated by constructing
the matrix \cite{hua1990matrix}
\begin{equation}
\tilde{L}=\Sigma_{\bar{R}}^{-1}U_{\bar{R}}^{*}YV_{\bar{R}}
\end{equation}
and computing its eigenvalues.

Indeed, if $\tilde{L}=\tilde{W}\bar{\Lambda}^{\#}\tilde{W}^{-1}$
where $\tilde{W},\bar{\Lambda}^{\#}\in\mathbb{C}^{\bar{R}\times\bar{R}}$
and $\bar{\Lambda}^{\#}$ is diagonal, the matrix pencil $\left(X,Y\right)$
can be approximated by
\begin{equation}
Y-\lambda X\approx U_{\bar{R}}\Sigma_{\bar{R}}\left(\tilde{L}-\lambda I\right)V_{\bar{R}}^{*}.
\end{equation}
Therefore, the columns of $V_{\bar{R}}\tilde{W}$ are its right generalized
eigenvectors and the rows of $\tilde{W}^{-1}\Sigma_{\bar{R}}^{-1}U_{\bar{R}}^{*}$
are its left generalized eigenvectors.

When comparing to the DMD matrix in (\ref{eq:DMD_matrix}) and (\ref{eq:eigendecomposition}),
we conclude that $\tilde{L}=\Sigma_{\bar{R}}\tilde{K}\Sigma_{\bar{R}}^{-1}$,
$\bar{\Lambda}^{\#}=\bar{\Lambda}$ and $\tilde{W}=\Sigma_{\bar{R}}^{-1}\tilde{V}$.
In other words, DMD and Matrix Pencil produce the same estimates for
eigenvalues and the reduced order matrices are similar (with $\Sigma_{\bar{R}}$
being the similarity transformation).

Finally we find that the estimates of both the adjoint DMD modes and
the left generalized eigenvectors of $\left(X,Y\right)$ are the same,
\begin{equation}
\tilde{V}^{-1}U_{\bar{R}}^{*}=\tilde{W}^{-1}\Sigma_{\bar{R}}^{-1}U_{\bar{R}}^{*},
\end{equation}
up to scaling and ordering.

Among closely related approaches are the State Space Method \cite{kung1983state}
and ESPRIT \cite{roy1986estimation} which construct the same Hankel
matrices and are equivalent to first order in the noise terms to the
Matrix Pencil method \cite{hua1991svd}. Other related approaches
are the SVD based Prony \cite{kumaresan1984prony} and the Prony Koopman
Mode Decomposition \cite{susuki2015prony} methods, which also yield
the same frequencies as DMD or the Matrix Pencil method for purely
sinusoidal data.

So far, our treatment of the $X$ and $Y$ data matrices was asymmetrical.
In our least squares approach we implicitly minimized the noise contribution
of the $X$ matrix alone when computing $\tilde{L}=\tilde{Y}\tilde{X}^{+}$.
However, Total Least Squares (TLS) variants of the above methods (TLS
Matrix Pencil \cite{hua1990total}, TLS ESPRIT \cite{van1991total},
TLS DMD \cite{hemati2017biasing}), have a symmetrical and
mutually similar treatment of the errors in both $X$ and~$Y$ data
matrices.

Another asymmetry arises due to our treatment of the data series as
going forward in time, even though this direction should not be preferred
over backward time propagation at least when the eigenvalues are all
on the unit circle. An interested reader may refer to the forward-and-backward
(FB) Matrix Pencil Method \cite{hua1990matrix} or the FB DMD \cite{dawson2016characterizing}
for an unbiased treatment of the data, although it is worth noticing
that the two methods are not the same. FB Matrix Pencil extends the
observable to include data which propagates both forward and backwards
in time, which is appropriate for sequential time series. On the other
hand, FB DMD first computes two propagation matrices, in forward and
backward time, and then computes the geometric average of the former
with the inverse of the later. At any rate, the TLS and FB versions
of the above algorithms all require at least twice the computational
effort of standard variants.

\subsection{\label{sub:ESPRIT}ESPRIT}
ESPRIT was originally developed for estimating the direction of arrival of signals generated by a set of narrowband emitters and collected by an array of sensors \cite{roy1986estimation}. As mentioned earlier, it is equivalent to the Matrix Pencil Method and  DMD with delayed observables when applied to the same sequential
data. However, an observation was made by the authors of ESPRIT, that when the emitted signals are uncorrelated their power can be estimated. This observation provided an inspiration for us to extend DMD to a clustering procedure for multiple times series
which will be described in the next section.

In the ESPRIT approach the data matrices are interpreted as
\begin{equation}
\begin{aligned}
X&=\underset{j=1}{\overset{R}{\sum}}a_{j}b_{j}^{T}=AB^{T}\\
Y&=\underset{j=1}{\overset{R}{\sum}}\lambda_{j}a_{j}b_{j}^{T}=A\Lambda B^{T}
\end{aligned},\label{eq:ESPRIT}
\end{equation}
where the columns of $X$ are outputs of a first array of sensors, and
the columns of $Y$ are outputs of a second array---identical, but
shifted relative to the first. Also, the $j$th column of $A=\begin{bmatrix}
a_{1} & \cdots & a_{R}\end{bmatrix}$ corresponds to direction-of-arrival vectors of (or sensors response
to) the $j$th signal, the $j$th columns of $B=\begin{bmatrix}
b_{1} & \cdots & b_{R}\end{bmatrix}$ consists of samples of the complex envelopes of the $j$th signal,
and $\Lambda=\mathrm{diag}\left\{ \lambda_{j}\right\} \in\mathbb{C}^{R\times R}$
corresponds to a phase shift between the first array of sensors and
the second one \cite{zoltowski1989sensor}.

Although presented in a different setup, those definition of the data
matrices are equivalent to the DMD matrices when constructed from
delayed observables. Indeed, from (\ref{eq:signal}), (\ref{eq:delayed_observable})
and (\ref{eq:data_matrices}) we see that
\begin{equation}
X=\underset{j=1}{\overset{R}{\sum}}\begin{bmatrix}
v_{j}\lambda_{j}^{t_{1}} & \cdots & v_{j}\lambda_{j}^{t_{M}}\\
v_{j}\lambda_{j}^{t_{1}+1} & \cdots & v_{j}\lambda_{j}^{t_{M}+1}\\
\vdots &  & \vdots\\
v_{j}\lambda_{j}^{t_{1}+d} & \cdots & v_{j}\lambda_{j}^{t_{M}+d}
\end{bmatrix},\:Y=\underset{j=1}{\overset{R}{\sum}}\begin{bmatrix}
v_{j}\lambda_{j}^{t_{1}+1} & \cdots & v_{j}\lambda_{j}^{t_{M}+1}\\
v_{j}\lambda_{j}^{t_{1}+2} & \cdots & v_{j}\lambda_{j}^{t_{M}+2}\\
\vdots &  & \vdots\\
v_{j}\lambda_{j}^{t_{1}+d+1} & \cdots & v_{j}\lambda_{j}^{t_{M}+d+1}
\end{bmatrix},
\end{equation}
and hence may choose
\begin{equation}
A=\begin{bmatrix}
\frac{v_{1}}{\left\Vert v_{1}\right\Vert }\lambda_{1}^{0} & \cdots & \frac{v_{R}}{\left\Vert v_{R}\right\Vert }\lambda_{R}^{0}\\
\frac{v_{1}}{\left\Vert v_{1}\right\Vert }\lambda_{1}^{1} & \cdots & \frac{v_{R}}{\left\Vert v_{R}\right\Vert }\lambda_{R}^{1}\\
\vdots &  & \vdots\\
\frac{v_{1}}{\left\Vert v_{1}\right\Vert }\lambda_{1}^{d} & \cdots & \frac{v_{R}}{\left\Vert v_{R}\right\Vert }\lambda_{R}^{d}
\end{bmatrix},\:B=\begin{bmatrix}
\left\Vert v_{1}\right\Vert \lambda_{1}^{t_{1}} & \cdots & \left\Vert v_{1}\right\Vert \lambda_{1}^{t_{M}}\\
\vdots &  & \vdots\\
\left\Vert v_{R}\right\Vert \lambda_{R}^{t_{1}} & \cdots & \left\Vert v_{R}\right\Vert \lambda_{R}^{t_{M}}
\end{bmatrix}\label{eq:ESPRIT_matrices}
\end{equation}
to put the matrices in the form of (\ref{eq:ESPRIT}).

In \cite{roy1986estimation} the authors show that when the sources
are uncorrelated their ``power'' can be estimated. In other words
one can estimate $B^{*}B$ if it is diagonal.  Unfortunately, while in the ESPRIT framework it is possible for the sources to be uncorrelated, in the formulation consistent with DMD (Eq.~(\ref{eq:ESPRIT_matrices})), $B^{*}B$ can never be diagonal because each data snapshot is always affected by all of the frequencies.
However,
in the next section we will consider data consisting of multiple signals,
each exhibiting a subset of frequencies with different magnitudes,
i.e.,
\begin{equation}
B=\begin{bmatrix}
b_{1,t_{1}}\lambda_{1}^{t_{1}} & \cdots & b_{1,t_{M}}\lambda_{1}^{t_{M}}\\
\vdots &  & \vdots\\
b_{R,t_{1}}\lambda_{R}^{t_{1}} & \cdots & b_{R,t_{M}}\lambda_{R}^{t_{M}}
\end{bmatrix}],
\end{equation}
where $b_{j,t_{k}}\in\mathbb{R}^{+}$. In that case DMD (or Matrix
Pencil method) might find the coefficients in $B$ and thus characterize
each time snapshot according to the prevalent frequencies it exhibits.
Furthermore, it might be possible to group ``similar'' snapshots
together as will be discussed next.

\section{\label{sec:Clustering}Clustering of time Series}

In this section we consider the prospects of DMD and Matrix Pencil
analysis to unsupervised classification of sequential data. Specifically,
we present a novel method for computing features of data series to
allow clustering them into groups of similar underlying dynamics,
that is, similar values of~$\lambda_j$ in (\ref{eq:signal}).

In section~\ref{sec:DMD_MP_comparison} the underlying assumption
was that all of the data comes from a single system, or a group of
independent systems to which DMD or Matrix Pencil methods are to be
applied separately. If we suspect that the set of time series constituting
the data comes from a smaller set of possible dynamical systems, we
may wish to group the series based on similar underlying dynamics.

Among the many approaches to clustering of time series, we focus on
ones that assume a model behind the data \cite{liao2005clustering}.
As an example, fitting autoregression parameters, defining some associated
metric and using them as features is one common approach \cite{piccolo1990distance,xiong2002mixtures}.
However, to the best of our knowledge, all contemporary clustering
techniques compute the model parameters (or features) for each time
series separately. Inspired by the ESPRIT algorithm (see discussion
in \ref{sub:ESPRIT}), we introduce a model-based approach which considers
all the times series at once, and uses elements of the DMD modes (or
generalized eigenvectors) as features for clustering purposes.

\subsection{\label{sub:DMD_clustering}DMD Clustering of Time Series}

Consider $N$ series each consisting of $d+2$ data points, $\left\{ y_{i}(t)\right\}\subset\mathbb{C}^{n}\:,i=1,\ldots,N,\:t=0,\ldots,d+1$.
We assume that each time series is a sum of damped sinusoids as in
section \ref{sec:DMD_MP_comparison} (Eq. (\ref{eq:signal})).  Each time
series was generated by one of several different dynamical systems, and the
number of different systems (which we call $P$) is significantly smaller than the
number of time series~$N$. Formally, we assume that the time series are split
into $P$ partitions
\begin{equation}
\begin{aligned}
I_{1}\cup\cdots\cup I_{P}&=\left\{ 1,\cdots,N\right\} \\
I_{1}\cap\cdots\cap I_{P}&=\emptyset
\end{aligned},\label{eq:partitions}
\end{equation}
each corresponding to the dynamics
\begin{equation}
x_{i}(t)=\underset{k=1}{\overset{l_{j}}{\sum}}v_{ijk}\lambda_{jk}^{t}\in\mathbb{C}^{n},\:i\in I_{j}\label{eq:systems}
\end{equation}
where $l_{j}$ is the number of oscillators in the $j$-th partition,
$\left\{ \lambda_{jk}\right\} _{k=1}^{l_{j}}\subset\mathbb{C}$ are
the corresponding complex frequencies for that partition and $x_{i}v_{ijk},\in\mathbb{C}^{n}$.

To state the general problem: we wish to estimate the complex frequencies
$\lambda_{jk}$ as well as the partitions $I_{j}$ from data given
by 
\begin{equation}
y_{i}(t)=x_{i}(t)+s_{i}(t),
\end{equation}
where $x_{i}(t)$ is given by (\ref{eq:systems}) and $s_{i}(t)$
is a noise term.

First, consider the left generalized eigenvectors of the matrix pencil
$\left(X,Y\right)$, where $X$ and $Y$ are data matrices
\begin{equation}
X=\begin{bmatrix}
y_{1}(0) & \cdots & y_{N}(0)\\
\vdots &  & \vdots\\
y_{1}(d) & \cdots & y_{N}(d)
\end{bmatrix},\:Y=\begin{bmatrix}
y_{1}(1) & \cdots & y_{N}(1)\\
\vdots &  & \vdots\\
y_{1}(d+1) & \cdots & y_{N}(d+1)
\end{bmatrix},\label{eq:data_multiple}
\end{equation}
of the dimensions $nd\times N$, which is the standard arrangement for DMD with delayed observables (similarly to (\ref{eq:sequential_data_matrices})).
The following property of the matrix
pencil in the noise-free case will be useful when determining the
partitions $I_{j}$. It requires that the number of series in each partition, $\big|I_j\big|$ (here, $|I_j|$ denotes the cardinality of the set $I_j$), is greater than $nl_{j}$, the dimensionality of the dynamic modes of the $j$th system.

\emph{Proposition 1}: For $j\in\{1,\ldots,P\}$, let $\Xi_j$ be a matrix with
$\big|I_j\big|$ columns given by
\[
\Xi_{j}=\begin{bmatrix}
v_{ij1}\\
\vdots\\
v_{ijl_{j}}
\end{bmatrix}_{i\in I_{j}}\in\mathbb{C}^{nl_{j}\times\left|I_{j}\right|},
\qquad i\in I_j.
\]
The matrix $\Xi_j$
thus has $nl_j$ rows and consists of all the coefficients (modes) of all systems
in (\ref{eq:systems}) belonging to the $j$th partition and stacked
\emph{on top of each other}. Assume that the matrices~$\Xi_{j}$ have full column
rank for all $j\in\{1,\dots,P\}$ (and in particular $\left|I_{j}\right|\ge nl_{j}=\mathrm{rank}\left(\Xi_{j}\right),\:\forall j$),
that $d\ge\sum l_{j}$ and that no noise is present, $s_{i}(t)=0,\:\forall i,t$.
Then the following hold:
\begin{enumerate}
\item The matrix pencil $\left(X,Y\right)$ has $nl$ linearly independent
left generalized eigenvectors where
\begin{equation}
l=\left|\left\{ \lambda_{jk} : 1\le j\le P,\:1\le k\le l_{j}\right\} \right|\label{eq:unique_evals}
\end{equation}
is the number of distinct eigenvalues among all systems in all partitions
($l\le\sum l_{j}$).
\item Let $p_{\lambda}$ be a generalized right eigenvector of the matrix
pencil $\left(X,Y\right)$ corresponding to the eigenvalue $\lambda$.
If the $i$th data snapshot corresponding to the $j$th partition
($i\in I_{j}$) doesn't exhibit $\lambda$ in its dynamics ($\lambda\notin\left\{ \lambda_{jk}\right\} _{k=1}^{l_{j}}$),
then the $i$th element of $p_{\lambda}$ is zero. 
\end{enumerate}
\emph{Proof}: See Appendix A.

In (\ref{eq:data_multiple}) one views the data as snapshots taken
at different times, similar to the approach in section~\ref{sec:DMD_MP_comparison}
(Eq.~(\ref{eq:sequential_data_matrices})). Alternatively, one can
arrange the data as
\begin{equation}
\hat{X}=\begin{bmatrix}
y_{1}(0) & \cdots & y_{1}(d)\\
\vdots &  & \vdots\\
y_{N}(0) & \cdots & y_{N}(d)
\end{bmatrix},\:\hat{Y}=\begin{bmatrix}
y_{1}(1) & \cdots & y_{1}(d+1)\\
\vdots &  & \vdots\\
y_{N}(1) & \cdots & y_{N}(d+1)
\end{bmatrix}\label{eq:data_mult_trans}
\end{equation}
with $\hat{X},\hat{Y}\in\mathbb{C}^{nN\times d}$, viewing all the
time series as consequent time snapshots of one large system. The
matrix pencil $\left(\hat{X},\hat{Y}\right)$ has similar properties
to $\left(X,Y\right)$ but requires ``less'' data and is therefore more useful for a numerical algorithm as stated in the next proposition and its discussion.

\emph{Proposition 2}: For $j\in \{1,\ldots,P\}$, let $\Xi_j$ be a matrix with
$n|I_j|$ rows, given by
\[
\Xi_j = \begin{bmatrix}
  v_{ij1} & \cdots & v_{ijl_{j}}
\end{bmatrix},
\qquad i\in I_{j}.\in\mathbb{C}^{n\left|I_{j}\right|\times l_{j}}
\]
Thus, the matrix $\Xi_j$ has $l_j$ columns and consists of all the coefficients (modes) of all systems
in (\ref{eq:systems}) belonging to the $j$th partition and placed
\emph{next to each other}. Assume that $\Xi_{j}$ has full row rank
for each $j$ (i.e. $n\left|I_{j}\right|\ge l_{j}=\mathrm{rank}\left(\Xi_{j}\right)$),
that $d\ge\sum l_{j}$ and that no noise is present, $s_{i}(t)=0,\:\forall i,t$.
Then the following hold:
\begin{enumerate}
\item The matrix pencil $\left(\hat{X},\hat{Y}\right)$ has $l$ linearly
independent left generalized eigenvectors where $l\le\sum l_{j}$
is the number of unique eigenvalues (see (\ref{eq:unique_evals})). 
\item Let $q_{\lambda}$ be a generalized left eigenvector of the matrix
pencil $\left(\hat{X},\hat{Y}\right)$ corresponding to the eigenvalue
$\lambda$. If the $i$th data snapshot corresponding to the $j$th partition
($i\in I_{j}$) doesn't exhibit $\lambda$ in its dynamics ($\lambda\notin\left\{ \lambda_{jk}\right\} _{k=1}^{l_{j}}$),
the $n$ elements of $q_{\lambda}$ beginning at $i\cdot n$ are all
zero. 
\end{enumerate}
\emph{Proof}: This follows directly from Proposition 1 for the matrix
pencil $\left(\hat{X}^{T},\hat{Y}^{T}\right)$ when considering each
$n$ dimensional system as $n$ one dimensional systems.

We note that the assumptions in Proposition 1 require more snapshots
and yield more generalized eigenvectors compared to Proposition 2.
For this reason we chose the arrangement in (\ref{eq:data_mult_trans})
over (\ref{eq:data_multiple}) in the example presented in section
\ref{sub:clustering_example}.

The presence of zero elements in the generalized eigenvectors
corresponding to different data series can be exploited in order to find
the partitions in (\ref{eq:partitions}). In other words, we propose
using the elements of the generalized eigenvectors (or adjoint DMD
modes) as features for clustering the data series.

Unfortunately, in presence of noise we do not expect any of the elements
of the generalized eigenvectors to be exactly zero as stated in Propositions
1 and~2. However, when using truncated SVD with the Matrix Pencil
method, the resulting generalized eigenvectors are continuous functions
of the data \cite{hua1990matrix}. Consequently, we expect the elements
which should theoretically be zero, to remain close to zero for
low levels of noise; this will be verified in an example in
section~\ref{sub:toy_example}.  We first discuss some details of
the implementation.

\subsection{\label{sub:cluster_impl}Implementation}

As mentioned earlier, the choice of arranging the data as in (\ref{eq:data_mult_trans})
has the benefits of giving a lower dimensional system and requiring
less data. We therefore will use the matrix pencil $\left(\hat{X},\hat{Y}\right)$
and its generalized right eigenvectors as features for clustering
the given time series.

\emph{Step 1}: Compute the Singular Value Decomposition (SVD): 
\begin{equation}
\hat{X}=U\Sigma V^{*}
\end{equation}
and choose a truncation value $\bar{l}$ for the singular values.
This gives $\hat{X}\approx U_{\bar{l}}\Sigma_{\bar{l}}V_{\bar{l}}^{*}$
where $U_{\bar{l}}\in\mathbb{C}^{nN\times\bar{l}}$, $\Sigma_{\bar{l}}\in\mathbb{C}^{\bar{l}\times\bar{l}}$
and $V\in\mathbb{C}^{d\times\bar{l}}$.

\emph{Step 2}: Construct $\hat{K}=U_{\bar{l}}^{*}\hat{Y}V_{\bar{l}}\Sigma_{\bar{l}}^{-1}\in\mathbb{C}^{\bar{l}\times\bar{l}}$
and find its left eigenvectors $\hat{W}\hat{K}=\hat{\Lambda}\hat{W}$,
where $\hat{\Lambda},\hat{W}\in\mathbb{C}^{\bar{l}\times\bar{l}}$
and $\hat{\Lambda}$ is diagonal.

\emph{Step 3}: Construct 
\begin{equation}
Q=\hat{W}U_{\bar{l}}^{*}\in\mathbb{C}^{\bar{l}\times nN},\label{eq:DMD_features_matrix}
\end{equation}
whose rows are generalized right eigenvectors of the matrix pencil
$\big(\hat{X},\hat{Y}\big)$.

\emph{Step 4}: We define DMD features as sub-matrices $q_{i}\in\mathbb{C}^{\bar{l}\times n},\:1\le i\le N$
of $Q$, whose elements are defined as 
\begin{equation}
\left(q_{i}\right)_{j,k}=\left|\left(Q\right)_{k,in+j}\right|,\:1\le j\le\bar{l},\:1\le k\le n,
\end{equation}
and $|\cdot|$ stands for the norm which is applied element-wise.
Furthermore there are many ways to define a metric on the above features,
for example
\begin{equation}
d\left(q_{i},q_{j}\right)=\left\Vert q_{i}-q_{j}\right\Vert _{F},
\end{equation}
where $\left\Vert \cdot\right\Vert _{F}$ is the Frobenius norm.

\emph{Step }5: Invoke any of the standard clustering methods on the
$\left\{ q_{i}\right\} _{i=1}^{N}$ features which have a one-to-one
correspondence with the $N$ time series.

We note that the matrix of right eigenvectors may also be expressed
as 
\begin{equation}
Q=\hat{W}\Sigma_{\bar{l}}^{-1}V_{\bar{l}}^{*}\hat{X}^{*}.
\end{equation}
Consequently, each of the DMD features can be seen as a transformation
of the appropriate data snapshot,
\begin{equation}
q_{i}=\left|\hat{W}\Sigma_{\bar{l}}^{-1}V_{\bar{l}}^{*}\begin{bmatrix}
y_{i}^{*}(0)\\
\vdots\\
y_{i}^{*}(d)
\end{bmatrix}\right|,\label{eq:DMD_features}
\end{equation}
where $\left|\cdot\right|$ stands for applying absolute value to
each elements of the matrix.

Moreover, as suggested previously in \cite{surana2018koopman}, DMD
features can be used during an unsupervised learning stage and later
applied to never-seen-before data. In our case, the ``training''
phase would consist of finding $\hat{W}$, $\Sigma_{\bar{l}}$, $V_{\bar{l}}$
and the centers of clusters in the feature space, $\mathbb{R}_{+}^{\bar{l}\times n}$.
Then, new data $\left\{ y(t)\right\} _{t=0}^{d}$ can be mapped into
the feature space via 
\begin{equation}
q=\left|\hat{W}\Sigma_{\bar{l}}^{-1}V_{\bar{l}}^{*}\begin{bmatrix}
y^{*}(0)\\
\vdots\\
y^{*}(d)
\end{bmatrix}\right|,
\end{equation}
and assigned to a previously ``learned'' cluster.

\section{\label{sec:Examples}Numerical Examples}

\subsection{\label{sub:toy_example} A Toy Example - Multiple 1D Sinusoidal Signals}

For illustration purposes we consider a set of one dimensional signals consisting of one or two sinusoids:
\begin{equation}
\begin{aligned}
  x_{1,...,6}(t)&=\alpha_{1,...,6}\exp(i\omega_{A}t)\\
  x_{7,...,12}(t)&=\alpha_{7,...,12}\exp(i\omega_{B}t)\\
  x_{13,...,23}(t)&=\alpha_{13,...,23}\exp(i\omega_{C}t)+\beta_{13,...,23}\exp(i\omega_{D}t),
\end{aligned}
\label{eq:ex1_dynamics}
\end{equation}
where $\omega_{A}=1.0$, $\omega_{B}=1.7$, $\omega_{C}=0.8$, $\omega_{D}=1.5$ and the coefficients $\alpha_{i},\beta_{i}\in\mathbb{C}$ are uniformly distributed on $\left\{ \left.1\le|z|\le2\right|z\in\mathbb{C}\right\}$.
The measurements $y(t)=\mathrm{Re}\left\{ x(t)\right\} +s(t)$ were taken at $t\in\left\{ 0,1,...,19\right\} $
and the noise term $s$ was normally distributed with $\sigma=0.1$. 

We now apply the technique presented in Section~\ref{sec:Clustering} to recover the frequencies from the data generated by \eqref{eq:ex1_dynamics} and cluster the signals based on those frequencies.
We arrange the data as in (\ref{eq:data_mult_trans}) into $\hat{X},\hat{Y}\in\mathbb{R}^{23\times19}$,
where $N=23$ is the number of signals.

\begin{figure}
\includegraphics[scale=0.9]{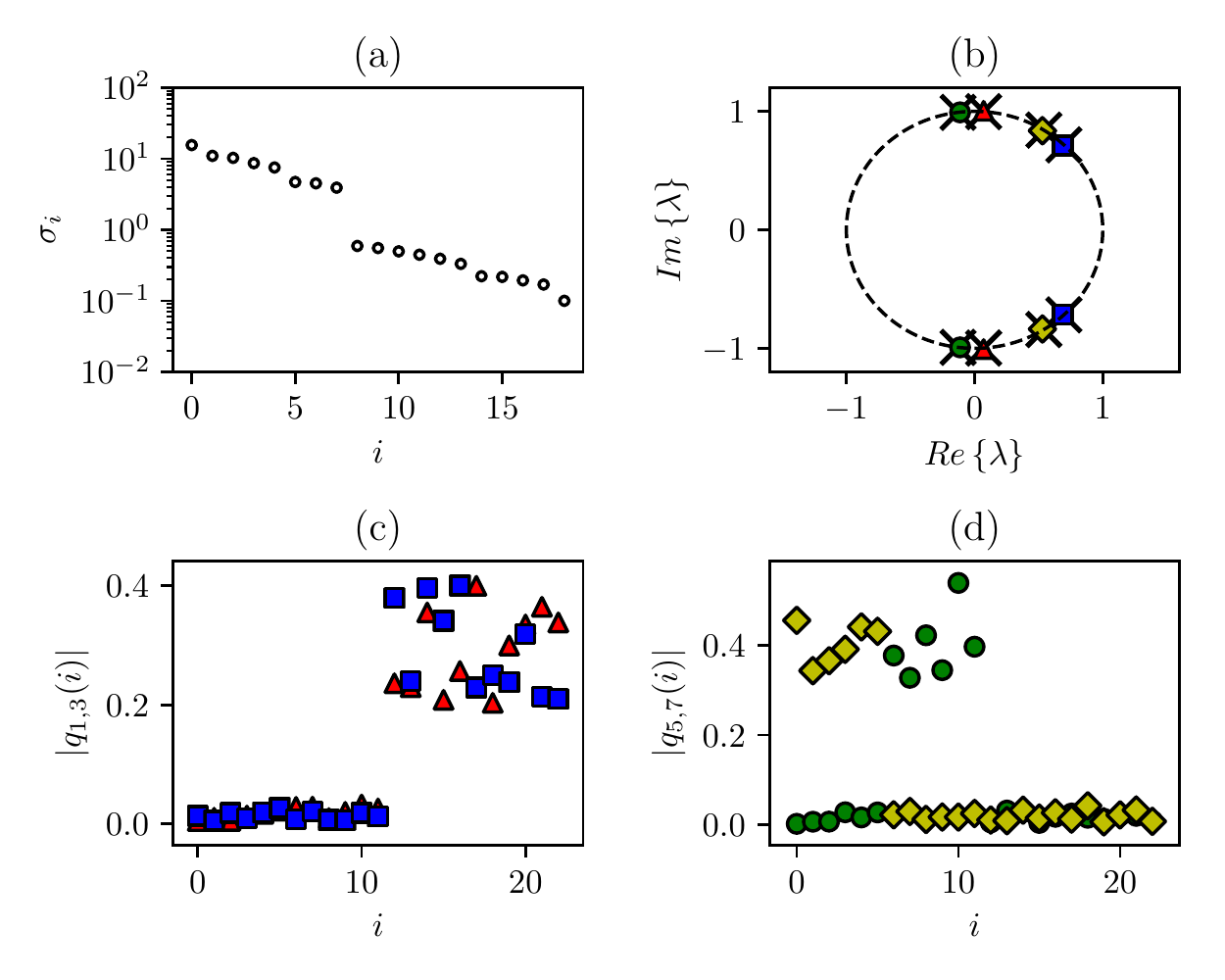}
 
\caption{(a) Singular values of the data matrix $\hat{X}$ generated by the dynamics in (\ref{eq:ex1_dynamics}) over 20 time steps. A drop in magnitude is clearly visible after the first eight, suggesting that the data contains that many complex frequencies (in this case, 4 complex conjugate pairs).
(b) Eigenvalues computed by DMD (colored shapes) as
approximations of the true eigenvalues ($\times$) of the discretization of (\ref{eq:ex1_dynamics}).
(c) Magnitudes of the first two complex conjugate pairs of generalized right eigenvectors of the data matrix pencil. Elements 1 to 12 are very close to zero indicating that the other sequences (13 to 23), exhibit the frequencies corresponding to these eigenvectors ($\omega_{C},\omega_{D}$).
(d) Magnitudes of the second and third complex conjugate  pairs of generalized right eigenvectors of the data matrix pencil. The frequency $\omega_{A}$ corresponds to the eigenvector with non-negligible elements at positions 1 through 6($\diamond$), which belong to the first six signals in the data. Similarly, $\omega_{B}$ corresponds to the eigenvector ($\circ$) with non-negligible elements belonging to signals 7 through 12.
}
\label{fig:example1}
\end{figure}

The singular values of the $\hat{X}$ are shown in figure~\ref{fig:example1}(a) and exhibit a sharp drop in magnitude after the first eight. This implies that the data contains four superimposed sinusoids. To estimate the frequencies we compute the DMD as suggested in section~\ref{sub:cluster_impl} with $\bar{l}=8$ modes. 

Figure~\ref{fig:example1}(b) shows the estimated eigenvalues of the matrix pencil $\big(\hat{X},\hat{Y}\big)$ which match the frequencies corresponding to the discrete measurements -- $\exp(\pm i\omega_{A,B,C,D})$. However, the singular values and the eigenvalues alone are not sufficient to determine which signals correspond to which frequencies. To do that, we compute the generalized right eigenvectors per (\ref{eq:DMD_features_matrix})

Figure~\ref{fig:example1}(c) and (d) show magnitudes  of four eigenvectors. Two of them (fig.~\ref{fig:example1}(c)) correspond to frequencies close to $\omega_{C}$ and $\omega_{D}$ and have elements with magnitudes close to zero at indices $1,...,12$. With proposition 2 in mind, this suggests that the other sequences ($13,...,23$) must exhibit those frequencies, which is indeed the case for the data generated by (\ref{eq:ex1_dynamics}). Similarly, the other two eigenvectors (fig.~\ref{fig:example1}(d)), corresponding to frequencies close to $\omega_{A}$ and $\omega_{B}$, have elements close to zero at all indices except $1,...,6$ and $7,...,12$ respectively. Again, we have correctly identified the signals those frequencies belong to. 

In this example one could easily conclude that a single frequency close to $\omega_{A}$ is present in sequences $1,...6$, $\omega_{B}$ in $7,...12$ and exactly two frequencies, $\omega_{C},\omega_{D}$, are present in sequences $13,...23$. For a large multidimensional dataset, this identification task requires a clustering algorithm as will be illustrated in the next example.

\subsection{\label{sub:clustering_example}Clustering of Lattice Regions in TEM
Image}

\begin{figure}
\includegraphics[scale=0.53]{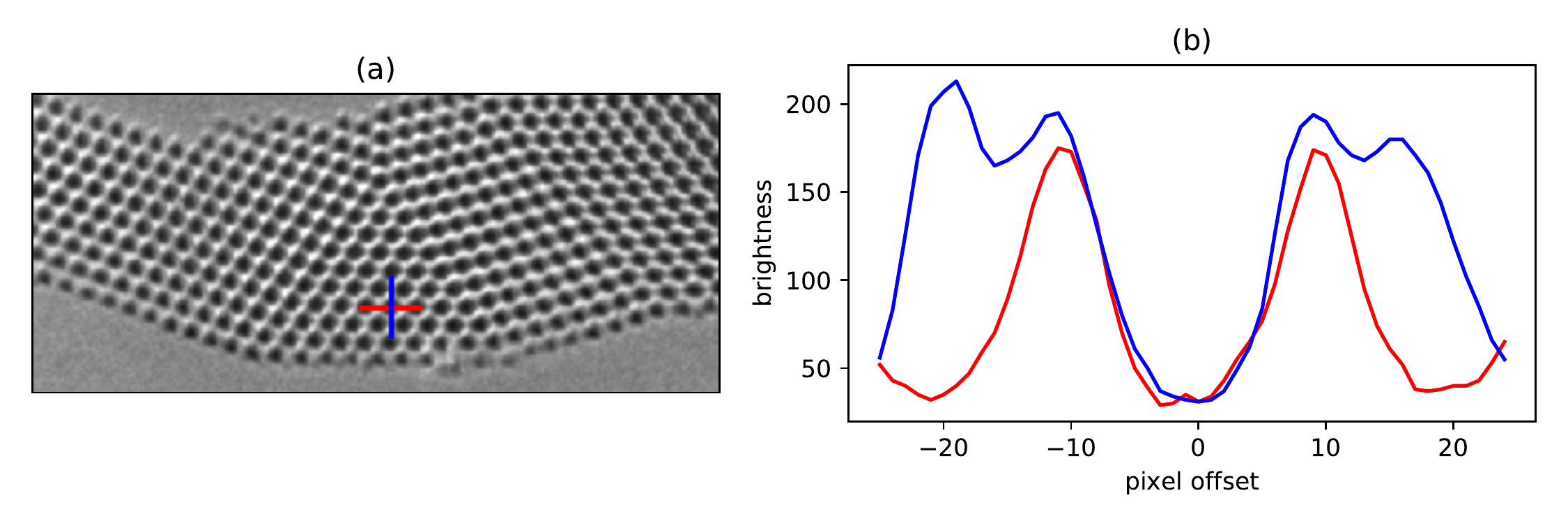}

\caption{\label{fig:lattice}Left: Transmission Electron Microscopy (TEM) lattice
images of Au atoms \cite{dahmen2009background} with scans along the
x and y axes around one pixel marked in red and blue respectively.
Right: the brightness profiles along the scans marked on the left
image. In the x direction (red), the period is about 19 pixels the
closest neighbor distance in the lattice. In the y direction (blue)
the period is about 28 pixels and it spans two rows.}
\end{figure}

For a more practical example, we consider a Transmission Electron Microscopy
(TEM) image of a lattice of gold atoms from \cite{dahmen2009background}, shown
in figure~\ref{fig:lattice}(a). We consider the variation in brightness
of the image as one transverses its pixels from left to right (x axis)
or top to bottom (y axis). Figure~\ref{fig:lattice}(b) shows those
variations over 51 pixels centered about a single pixel in a particular
region of the lattice. In that region, the hexagonal close-packed
(HCP) lattice is oriented (almost) parallel to the x axis of the image.
The brightness variation in the x direction has a period of about
19 pixels - the closest neighbor distance in the lattice. The variation
in the y axis, in this region, has a period of exactly $1.5$ times
the period in the x axis (about 28 pixels). One period in the y axis
however includes \emph{two} rows of the lattice, hence the variations
in the y direction in that region exhibit two frequencies,
$\omega_{y}=1.5\omega_{x}\approx 1/28$
and $2\omega_{y}$.

In various regions in figure~\ref{fig:lattice}(a), the grain has
different orientations and therefore different frequencies of variations
in brightness in x and y directions. In this example will use DMD features and constrained hierarchical clustering \cite{openshaw1977geographical}
to split the pixels of the image based on these variations.

We construct the data matrices by combining the brightness data
(figure~\ref{fig:lattice}(b)) for each pixel in each direction x
and y ($n=2$).  Letting $P_{ij}$ denote the brightness of pixel $(i,j)$, we
arrange the data as
\begin{equation}
  \left\{ y_{ij}(t)\right\} _{t=0}^{d}=\left\{
    \begin{bmatrix}
      P_{i-d/2,j}\\
      P_{i,j-d/2}
    \end{bmatrix},...,
    \begin{bmatrix}
      P_{i-1,j}\\
      P_{i,j-1}
    \end{bmatrix},
    \begin{bmatrix}
      P_{i,j}\\
      P_{i,j}
    \end{bmatrix},
    \begin{bmatrix}
      P_{i+1,j}\\
      P_{i,j+1}
    \end{bmatrix},...,
    \begin{bmatrix}
      P_{i+d/2,j}\\
      P_{i,j+d/2}
    \end{bmatrix}\right\} 
\end{equation}
The number of data series ($N$) is the number of pixels for which
brightness variations over 51 pixels ($d=50$) were collected. The
partitions in (\ref{eq:partitions}) therefore correspond to regions
where grain orientation remains the same or regions where no lattice is present.

\begin{figure}
\includegraphics[scale=0.59]{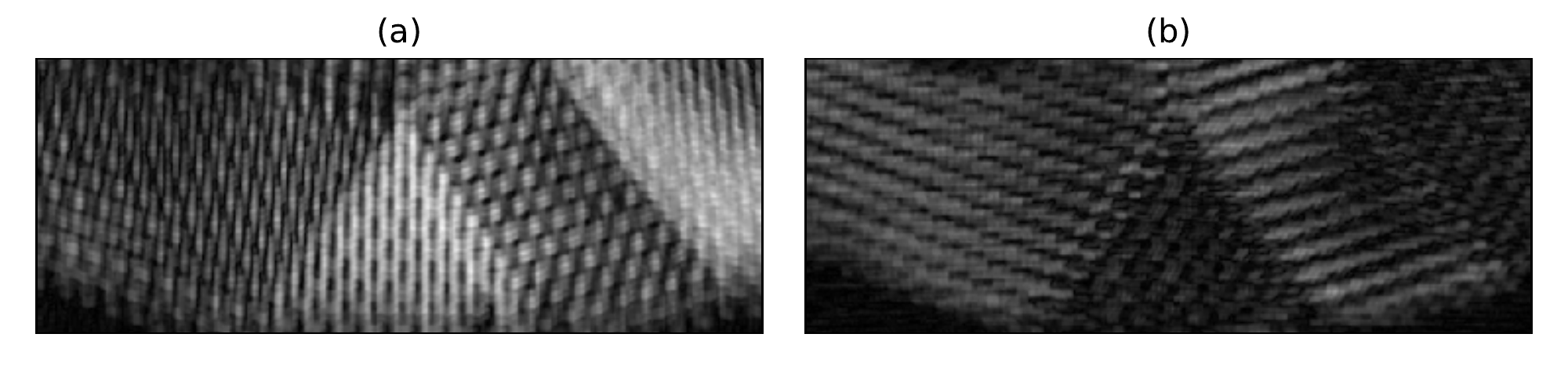}

\caption{\label{fig:lattice_mode} The adjoint DMD mode (or generalized left
eigenvector) based on data shown in figure \ref{fig:lattice}(a) and
corresponding to the frequency $0.067\approx2\omega_{y}\approx\frac{1}{14}$.
The mode was split into elements corresponding to the x (right) and
y (left) scanning directions. The elements in the triangular shaped region in the middle
have large magnitudes for the y direction (left) and are
close to zero for the x direction (right) as predicted by figure \ref{fig:lattice}(b).}
\end{figure}

The generalized left eigenvectors of $\big(\hat{X},\hat{Y}\big)$
in (\ref{eq:data_mult_trans}) have $2N$ elements (per \ref{sub:cluster_impl})
which we rearrange back into the shape of the original image, once
for the x and once for the y direction. Figure~\ref{fig:lattice_mode}
shows the element-wise absolute value of a generalized left eigenvector
corresponding to the frequency closest to $2\omega_{y}\approx 1/14$.
As discussed earlier, this frequency is present only in brightness
variations in the y direction in the triangular shaped region where
the lattice is parallel to the x axis (marked in figure~\ref{fig:lattice}(a)).
Indeed figure~\ref{fig:lattice_mode}(a) shows the part of
the eigenvector corresponding to the y direction, which exhibits large
magnitudes in the triangular area where this frequency occurs. On
the other hand, the elements corresponding to the $2\omega_{y}$ frequency
but the x direction (figure~\ref{fig:lattice_mode}(b)) are close
to zero in the same area as expected per Proposition 2 in
\ref{sub:DMD_clustering}.

Having computed the feature vectors (Eq.~(\ref{eq:DMD_features})),
we may proceed with a clustering scheme of our choice. Since in this
particular case the feature space alone doesn't contain any spatial
information (i.e., pixels relative position), we speculate
that a contiguity-constrained clustering \cite{murtagh1985survey}
approach would be appropriate. Specifically, a constrained version
\cite{openshaw1977geographical} of the Ward hierarchical agglomerative
clustering algorithm \cite{ward1963hierarchical} seems to give satisfactory
results. In this approach, each pixel begins as a singleton cluster.
Then clusters are iteratively merged in a greedy manner, such that
each two newly merged clusters minimize the Ward criterion \cite{ward1963hierarchical}
with respect to the DMD features. However, clusters can be merged only
if they have adjacent pixels, i.e., if they satisfy the spatial connectivity
constraints. This reduces computational time while ensuring that regions
of the image corresponding to similar dynamics remain connected.
An implementation of the Ward agglomerative clustering scheme is available
with \emph{Scikit-learn}~\cite{pedregosa2011scikit}, a Python machine
learning library.

\begin{figure}
\includegraphics[scale=0.77]{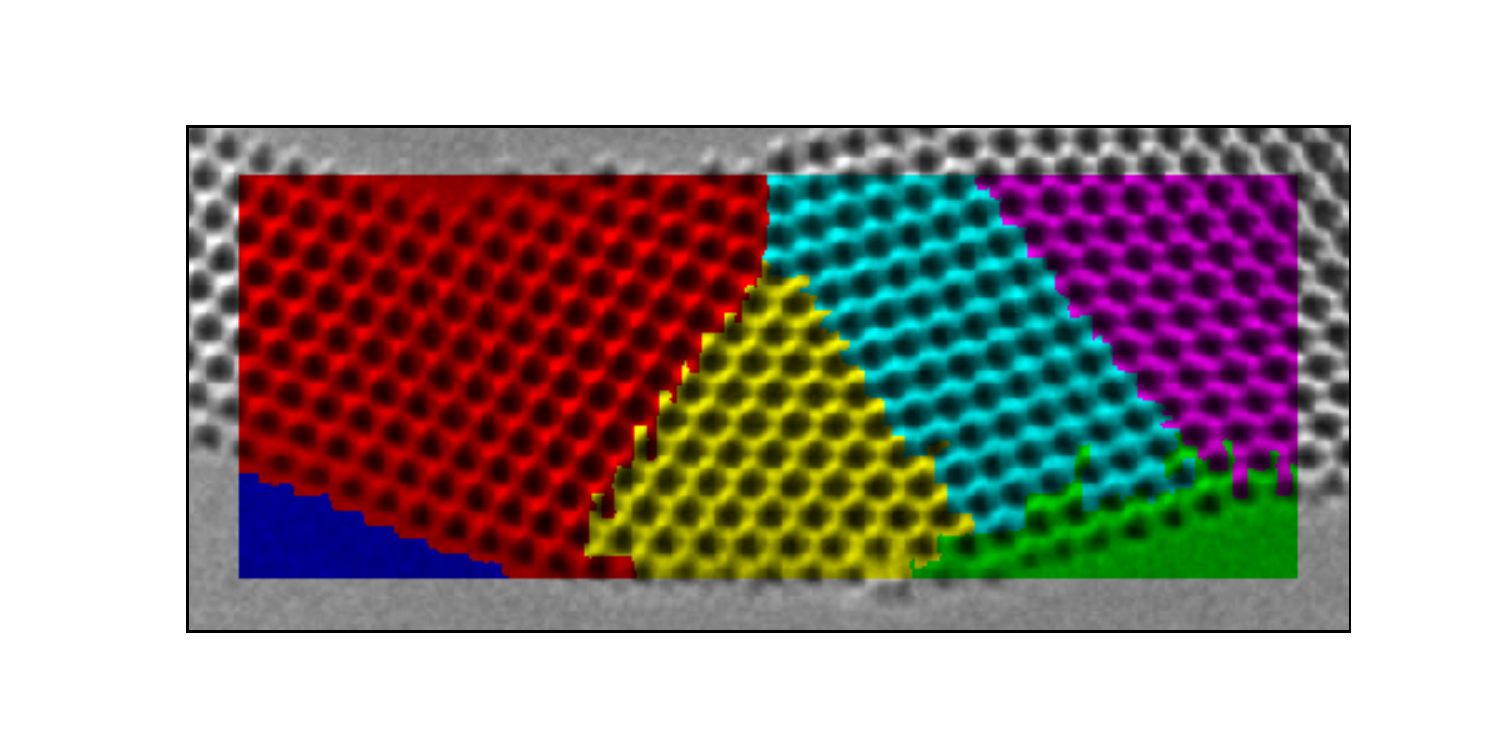}

\caption{\label{fig:clustering}Clustering of pixels based on DMD features
computed from the brightness variation in x and y direction (see figure \ref{fig:lattice}).
A choice of 6 clusters gives the correct separation into regions with
various lattice orientations and spatially separated regions with no lattice.}
\end{figure}

The results are shown on figure~\ref{fig:clustering}, where the pixels
are split into 6 regions based on their DMD features computed as specified
in section~\ref{sub:cluster_impl}. Since hierarchical clustering produces
a dendogram which has separate pixels as leaves and a single cluster
containing the whole image as a root, any number of clusters between
those two extremes is a priori valid. However, once the dendogram is
complete, it requires little additional computational effort to retrieve
the partitions for any possible number of clusters. In the example
presented here, we find that a choice of 6 partitions seem to capture
the different regions in figure \ref{fig:lattice} (a) correctly.
Indeed in figure~\ref{fig:clustering} we see four partitions where the
lattice orientation causes different frequencies in variation in brightness
in x and y direction. Two more partitions have no lattice present
and form two spatially separated regions.

\section{Conclusions}

We have described connections between Dynamic Mode Decomposition (DMD) and a family of frequency detection
methods (namely, the Matrix Pencil, State Space, and ESPRIT methods) that were
formerly known to be equivalent to one another~\cite{hua1991svd}. Furthermore, it
was shown that, for sequential data, DMD with delayed observables
yields exactly the same eigenvalues as the Matrix Pencil method, while
the generalized eigenvectors were found to be the adjoint DMD modes.
Based on these similarities, and inspired by the formulation of the
ESPRIT method, a novel approach for clustering of time series was
proposed.

We have considered data in the form of a large number of time series
with a smaller number of possible underlying dynamics (frequencies)
for each series. It was shown that, when arranging all the data together,
the generalized eigenvectors (adjoint DMD modes) of the corresponding
Hankel matrices have an interesting property: their elements corresponding
to certain series and frequencies tend to have magnitudes close to
zero when those series do not exhibit those particular frequencies.
We therefore conclude that projecting a sequence onto DMD modes computed
from a larger set of given sequences, is a viable feature extraction
method given a large amount of unlabeled series.

While existing approaches to clustering of time series require extracting
features from each sequence first, DMD and the Matrix Pencil methods
estimate the number of required features and extract them across all
series at once. Various metrics can then be defined on the elements
of the adjoint DMD modes and a range of clustering algorithms can
be applied on the time series accordingly. This method has been illustrated
by clustering regions of distinct patterns in an image based on the
variance in brightness of spatially close pixels.

\section*{Appendix A}

\emph{Proof of Proposition 1}:

Without loss of generality, we assume throughout the proof that the
columns of $X$ and $Y$ defined in (\ref{eq:data_multiple}) are
grouped together based on the partitions, i.e. the first $\left|I_{1}\right|$
columns are data from systems in the first partition etc. This makes
the proof more tractable but doesn't affect our conclusion regarding
the elements of the generalized eigenvectors.

First, consider the case when $n=1$, $\left|I_{j}\right|=l_{j}$
and $\Xi_{j}=I_{l_{j}\times l_{j}}$ is the identity matrix. The data
matrices take the following form:
\begin{equation}
X^{\#}=\begin{bmatrix}
1 & 1 & \cdots & 1 & 1 & \cdots & 1\\
\lambda_{11} & \lambda_{12} & \cdots & \lambda_{1l_{1}} & \lambda_{21} & \cdots & \lambda_{Pl_{P}}\\
\vdots & \vdots &  & \vdots & \vdots &  & \vdots\\
\lambda_{11}^{d} & \lambda_{12}^{d} & \cdots & \lambda_{1l_{1}}^{d} & \lambda_{21}^{d} & \cdots & \lambda_{Pl_{P}}^{d}
\end{bmatrix}\in\mathbb{C}^{d\times\sum l_{j}},
\end{equation}
where the columns are increasing powers of all the complex frequencies
of all systems arranged in order of their partitions.

From the statement of the theorem $d\ge l$, where $l\le\sum l_{j}$
is the number of unique exponents among all partitions, thus the rank
of $X^{\#}$ is $l$. Let $e_{jk}$ be vectors of the standard basis
of $\mathbb{C}^{\sum l_{j}}$ such that the non-zero element of $e_{jk}$
corresponds to the column of $X^{\#}$ where the powers of $\lambda_{jk}$
appear, and consider the matrix pencil
\begin{equation}
\lambda X^{\#}-Y^{\#}=\begin{bmatrix}
\lambda-\lambda_{11} & \cdots & \lambda-\lambda_{1l_{1}} & \lambda-\lambda_{21} & \cdots\\
\vdots &  & \vdots & \vdots\\
\lambda_{11}^{d}\left(\lambda-\lambda_{11}\right) & \cdots & \lambda_{1l_{1}}^{d}\left(\lambda-\lambda_{1l_{1}}\right) & \lambda_{21}^{d}\left(\lambda-\lambda_{21}\right) & \cdots
\end{bmatrix},
\end{equation}
where $Y^{\#}$ is defined similarly to $X^{\#}$ with all powers increased
by one.

Clearly, if $\lambda$ is one of the exponents $\lambda_{jk}$, it
is a generalized eigenvalue of the matrix pencil $\left(X^{\#},Y^{\#}\right)$,
since in that case $\mathrm{rank}\left(Y^{\#}-\lambda X^{\#}\right)=\mathrm{rank}\left(X^{\#}\right)-1$.
The corresponding generalized right eigenvector is
\begin{equation}
p_{\lambda}^{\#}=\underset{\lambda_{jk}=\lambda}{\sum}e_{jk}
\end{equation}
since it can be easily shown that $p_{\lambda}^{\#}$ in the row space
of $X^{\#}$. We have therefore found all the $l$ generalized eigenvectors
of $\left(X^{\#},Y^{\#}\right)$. Moreover, generalized eigenvectors
that belong to distinct eigenvalues are orthogonal in this case (whether
they belong to the same partition or not).

Remaining in the one dimensional case ($n=1$), we now relax the assumptions
on the number of series and modes per partition, allowing $\left|I_{j}\right|\ge l_{j}$
and $\Xi_{j}$ of any form as long as satisfies $\mathrm{rank}\left(\Xi_{j}\right)=l_{j}$.
The data matrices may now be expressed as
\begin{equation}
X=\begin{bmatrix}
\begin{array}{ccc}
1 & \cdots & 1\\
\vdots &  & \vdots\\
\lambda_{11}^{d} & \cdots & \lambda_{1l_{1}}^{d}
\end{array} & \cdots & \begin{array}{ccc}
1 & \cdots & 1\\
\vdots &  & \vdots\\
\lambda_{P1}^{d} & \cdots & \lambda_{Pl_{P}}^{d}
\end{array}\end{bmatrix}\begin{bmatrix}
\Xi_{1} &  & 0\\
 & \ddots\\
0 &  & \Xi_{P}
\end{bmatrix},
\end{equation}
and the matrix pencil as
\begin{equation}
Y-\lambda X=\left(Y^{\#}-\lambda X^{\#}\right)\begin{bmatrix}
\Xi_{1} &  & 0\\
 & \ddots\\
0 &  & \Xi_{P}
\end{bmatrix}
\end{equation}

Since $\Xi_{j}\in\mathbb{C}^{l_{j}\times\left|I_{j}\right|},\;\forall j$
have full column rank,
\begin{equation}
p_{\lambda}=\begin{bmatrix}
\Xi_{1} &  & 0\\
 & \ddots\\
0 &  & \Xi_{P}
\end{bmatrix}^{+}p_{\lambda}^{\#}=\underset{\lambda_{jk}=\lambda}{\sum}\begin{bmatrix}
\Xi_{1}^{+} &  & 0\\
 & \ddots\\
0 &  & \Xi_{P}^{+}
\end{bmatrix}e_{jk}\label{eq:eig_structure}
\end{equation}
is clearly a generalized eigenvector of $\left(X,Y\right)$ when $\lambda$
is one of the $l$ unique exponents. Since $X$ has rank $l$, we
have found all the generalized right eigenvectors. Moreover the structure
of $e_{jk}$ and the matrix on the right hand side of (\ref{eq:eig_structure}),
reveals that $p_{\lambda}$ may only have non-zero elements at indices
belonging to partitions which exhibit $\lambda$ in their dynamics.
Thus, the proposition is proved for this case.

Finally, we consider the case $n>1$ in which $\Xi_{j}\in\mathbb{C}^{nl_{j}\times\left|I_{j}\right|}$
has rank $nl_{j}$. In this case, the data matrices can be expressed
as
\begin{equation}
X=\left(X^{\#}\otimes I_{n\times n}\right)\begin{bmatrix}
\Xi_{1} &  & 0\\
 & \ddots\\
0 &  & \Xi_{P}
\end{bmatrix},\:Y=\left(Y^{\#}\otimes I_{n\times n}\right)\begin{bmatrix}
\Xi_{1} &  & 0\\
 & \ddots\\
0 &  & \Xi_{P}
\end{bmatrix}
\end{equation}
where $\otimes$ denotes the Kronecker product and $I_{n\times n}$
and an $n$ by $n$ identity matrix.

The argument in the proof of the $n=1$ case applies to each dimension
separately, resulting in $n$ times as many generalized eigenvalues
and eigenvectors, and the conclusion still holds.

 \bibliographystyle{abbrv}
\bibliography{DMD_Clustering}

\end{document}